\documentclass[12pt,twoside,reqno]{amsart}
\usepackage{amsmath,times,epsfig,amssymb,amsbsy,amscd,amsfonts,amstext,graphicx}

\calclayout
\newtheorem{lem}{\bf Lemma}
\newtheorem{prop}{\bf Proposition}
\newtheorem{cor}{\bf Corollary}
\newtheorem{rk}{\bf Remark}

\newtheorem{defn}{\bf Definition}

\begin{document}
\title{Relative Garside elements of Artin Monoids}
\author{Usman Ali$^{1}$, Barbu Berceanu$^{2}$, Zaffar Iqbal$^{3}$}
\address{$^{1}$ Center for Advanced Studies in Pure and Applied Mathematics, BZU, Multan, Pakistan.}
\email{uali@bzu.edu.pk}
\address{$^{2}$ Abdus Salam School of Mathematical Sciences, GC University, Lahore,
Pakistan and Institute of Mathematics Simion Stoilow, Bucharest, Romania (permanent address).}
\email{Barbu.Berceanu@imar.ro}
\address{$^{3}$ Department of Mathematics, University of Gujrat, Pakistan}
\email{zaffarsms@yahoo.com}
\thanks{ This research is partially supported by Higher Education Commission, Pakistan.\\
2010 AMS classification: Primary 20F36; Secondary 20M05, 20F10.}

\begin{abstract}
We introduce a relative Garside element, the quotient of the corresponding Garside
elements $\Delta(\Gamma_{n-1})$ and $\Delta(\Gamma_{n})$, for a pair of Artin monoids 
associated to Coxeter graphs $\Gamma_{n-1}\subset\Gamma_{n}$, the second graph containing 
a new vertex. These relative elements give a recurrence relation between Garside elements. 
As an application, we compute explicitly the Garside elements of Artin monoids
corresponding to spherical Coxeter graphs or the longest elements of the associated 
finite Coxeter groups. 

\end{abstract}

\maketitle

\fontsize{11}{20}\selectfont

\pagestyle{myheadings} \markboth{\centerline {\scriptsize U. Ali,
B. Berceanu, Z. Iqbal}} {\centerline {\scriptsize Garside elements
of Artin monoids}}

\section{Introduction}\label{intro}

The Garside element for the braid group (and for the braid monoid) corresponding to Artin's
presentation (see~\cite{art},~\cite{mor})
$$\mathcal{B}_{n}\!=\!\left\langle x_{1},x_{2},\ldots,x_{n-1}\,\Bigg|
\begin{array}{l}
x_{i}x_{j}=x_{j}\,x_{i}\,\,\,\mbox{if} \, \mid i-j\mid\,\,\geq2 \\
x_{i+1}\,x_{i}\,x_{i+1}=x_{i}\,x_{i+1}\,x_{i}\mbox{ if } 1\leq i\leq n-2
\end{array}
\right\rangle $$
is given by (see~\cite{gar})
$$\Delta_{n}=x_{1}(x_{2}x_{1})(x_{3}x_{2}x_{1})\dots (x_{n-1}x_{n-2}\dots x_{1})$$
(to represent $\Delta_{n}$, or, more general, to represent an element of a monoid
as a product of generators, we chose the smallest word in the
length-lexicographic order induced by the order between generators
$x_1<x_2<\ldots <x_{n-1}$). For other representations of $\Delta_{n}$, including
the most used formula, see~\cite{BP}. The right
quotient $\Delta_{n}^{-1}\Delta_{n+1}=x_nx_{n-1}\ldots x_1$ will be called
\textit{the relative Garside element} corresponding to the embedding of the Coxeter
graphs $A_{n-1}\subset A_{n} $:
\begin{center}
\begin{picture}(360,46)
\thicklines                   \put(-5,25){$A_{n-1}:$}       \put(185,25){$A_{n}:$}
\put(330,27){\line(1,0){30}}  \put(358,17){$x_{n}$}          \put(358,24){$\bullet$}
\put(170,25){$\subset$}       \multiput(0,0)(180,0){2}{\put(30,27){\line(1,0){120}}
\put(28,24){$\bullet$}        \put(58,24){$\bullet$}         \put(88,24){$\bullet$}
\put(148,24){$\bullet$}       \put(28,17){$x_{1}$}           \put(58,17){$x_{2}$}
\put(88,17){$x_{3}$}          \put(144,17){$x_{n-1}$}        \put(116,16){$\dots$}}
\end{picture}
\end{center}
This element and its left divisors play a
central role in the construction of the Gr\"obner basis for the classical braid
monoid (see~\cite{tez},~\cite{bok}). In this paper we give few characterizations of the
relative Garside element corresponding to an extension of Coxeter graphs by one
new vertex, see Propositions 1- 4, and we will use these to compute inductively the
Garside elements of the Artin monoids of spherical type, see Corollaries 1-7.

We start to recollect some facts about Coxeter graphs, Coxeter and Artin groups,
Artin monoids, and Garside elements. Let $S$ be a set. A \emph{Coxeter matrix} over
$S$ is a square matrix $M=(m_{st})_{s,t\in S}$ indexed by the elements of $S$
such that

$m_{ss}=1$ for all $s\in S$ and $m_{st}=m_{ts}\in\{2,3,4,\ldots,\infty\}$ for all 
$s,t\in S,\,s\neq t$.

\noindent The associated \emph{Coxeter graph} $\Gamma=\Gamma(M)$ is a
labeled graph defined by the following data:

$S$ is a set of vertices of $\Gamma$;

two vertices $s,t\in S$ are joined by an edge if $m_{st}\geq3$,
with label $m_{st}$ if $m_{st}\geq4$.

\noindent A Coxeter matrix $M=(m_{st})_{s,t\in S}$ is usually represented by
its Coxeter graph $\Gamma(M)$.

\begin{defn}
Let $M=(m_{st})_{s,t\in S}$ be the Coxeter matrix of the Coxeter graph $\Gamma$. Then
the group defined by
$$\mathcal{W}(\Gamma)=\big\langle s\in S\mid (st)^{m_{st}}=1\mbox{ for all }
s,t\in S\mbox{ satisfying }m_{st}\neq\infty\big\rangle $$
is called the \emph{Coxeter group} of type $\Gamma$.
\end{defn}

In an equivalent way we can write $\mathcal{W}(\Gamma)=\big\langle s\in S\mid
s^{2}=1,\mathop{\underbrace{sts\dots}} \limits_{m_{st}\,\text{factors}}=
\mathop{\underbrace{tst\dots}} \limits_{m_{st}\,\text{factors}}\big\rangle$.

We call $\Gamma$ to be of \emph{spherical type} if $\mathcal{W}(\Gamma)$ is
a finite group. A graph is of spherical type if and only if it has finitely
many connected components, any of them from the next list (see~\cite{bbk}):

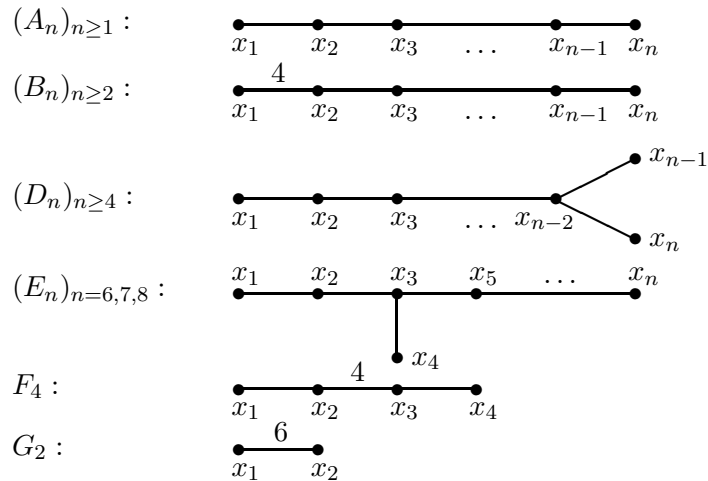
\begin{center}
\begin{picture}(300,46)
\thicklines              \put(100,27){\line(1,0){150}}     \put(15,25){$(A_{n})_{n\geq1}:$}
\put(98,24){$\bullet$}   \put(128,24){$\bullet$}           \put(158,24){$\bullet$}
\put(218,24){$\bullet$}  \put(248,24){$\bullet$}
\put(98,17){$x_{1}$}     \put(128,17){$x_{2}$}             \put(158,17){$x_{3}$}
\put(218,17){$x_{n-1}$}  \put(248,17){$x_{n}$}             \put(186,16){$\dots$}
\end{picture}

\begin{picture}(300,24)
\thicklines             \put(100,27){\line(1,0){150}}      \put(15,25){$(B_{n})_{n\geq2}:$}
\put(98,24){$\bullet$}  \put(128,24){$\bullet$}            \put(158,24){$\bullet$}
\put(218,24){$\bullet$} \put(248,24){$\bullet$}            \put(113,30){$4$}
\put(98,17){$x_{1}$}    \put(128,17){$x_{2}$}              \put(158,17){$x_{3}$}
\put(218,17){$x_{n-1}$} \put(248,17){$x_{n}$}              \put(186,16){$\dots$}
\end{picture}

\begin{picture}(300,40)
\thicklines            \put(100,27){\line(1,0){120}}       \put(220,26.5){\line(2,1){30}}
\put(98,24){$\bullet$} \put(128,24){$\bullet$}             \put(158,24){$\bullet$}
\put(218,24){$\bullet$}\put(248,39.5){$\bullet$}           \put(248,9){$\bullet$}
\put(98,17){$x_{1}$}   \put(128,17){$x_{2}$}               \put(158,17){$x_{3}$}
\put(205,17){$x_{n-2}$}\put(256,9){$x_{n}$}                \put(256,41){$x_{n-1}$}
\put(186,16){$\dots$} \put(15,25){$(D_{n})_{n\geq4}:$}    \put(220.5,26.5){\line(2,-1){30}}
\end{picture}

\begin{picture}(300,35)
\thicklines             \put(100,27){\line(1,0){150}}        \put(160.5,27){\line(0,-1){24}}
\put(98,24){$\bullet$}  \put(128,24){$\bullet$}              \put(158,24){$\bullet$}
\put(188,24){$\bullet$} \put(248,24){$\bullet$}              \put(158,0){$\bullet$}
\put(98,32){$x_{1}$}    \put(128,32){$x_{2}$}                \put(158,32){$x_{3}$}
\put(166,0){$x_{4}$}    \put(188,32){$x_{5}$}                \put(248,32){$x_{n}$}
\put(216,31){$\dots$}   \put(15,26){$(E_{n})_{n=6,7,8}:$}
\end{picture}

\begin{picture}(300,35)
\thicklines             \put(100,27){\line(1,0){90}}         \put(98,24){$\bullet$}
\put(128,24){$\bullet$} \put(158,24){$\bullet$}              \put(188,24){$\bullet$}
\put(98,17){$x_{1}$}    \put(128,17){$x_{2}$}                \put(158,17){$x_{3}$}
\put(188,17){$x_{4}$}   \put(143,30){$4$}                    \put(15,25){$F_{4}:$}
\end{picture}

\begin{picture}(300,22)
\thicklines             \put(100,27){\line(1,0){30}}         \put(98,24){$\bullet$}
\put(128,24){$\bullet$} \put(98,17){$x_{1}$}                 \put(128,17){$x_{2}$}
\put(114,30){$6$}       \put(15,25){$G_{2}:$}
\end{picture}

\begin{picture}(300,32)
\thicklines            \put(100,27){\line(1,0){60}}          \put(98,24){$\bullet$}
\put(128,24){$\bullet$}\put(158,24){$\bullet$}               \put(98,17){$x_{1}$}
\put(128,17){$x_{2}$}  \put(158,17){$x_{3}$}                 \put(113,30){$5$}
\put(198,24){$\bullet$}\put(228,24){$\bullet$}               \put(258,24){$\bullet$}
\put(288,24){$\bullet$}\put(198,17){$x_{1}$}                 \put(228,17){$x_{2}$}
\put(258,17){$x_{3}$}  \put(288,17){$x_{4}$}                 \put(178,17){$,$}
\put(213,30){$5$}      \put(15,25){$(H_{n})_{n=3,4}:$}       \put(200,27){\line(1,0){90}}
\end{picture}

\begin{picture}(300,32)
\thicklines              \put(100,37){\line(1,0){30}}
\put(98,34){$\bullet$} \put(128,34){$\bullet$}
\put(98,27){$x_{1}$}                \put(128,27){$x_{2}$}
\put(114,42){$p$}
\put(15,35){$\big(I_{2}(p)\big)_{p\geq5,p\neq6}:$}
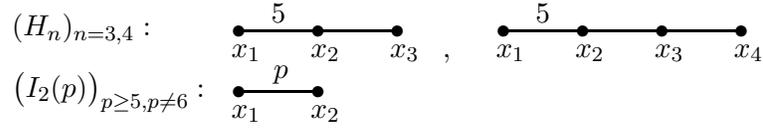
\put(15,0){\mbox{\textbf{Figure 1.1}. The connected spherical type
Coxeter graphs.}}
\end{picture}
\end{center}

\begin{defn}
If $\Gamma$ is a Coxeter graph, its associated \emph{Artin group} is defined by
$$\mathcal{A}(\Gamma)=\big\langle s\in S\,\big|\mathop{\underbrace{sts\dots}}
\limits_{m_{st}\,\text{factors}}=\mathop{\underbrace{tst\dots}}\limits_{m_{st}\,
\text{factors}}\mbox{ for all }s,t\in S\mbox{ satisfying }m_{st}\neq\infty\big\rangle.$$
\end{defn}

The set of positive elements (i.e. the elements which are product of generators
with positive exponents) in an Artin group $\mathcal{A}(\Gamma )$ is called the
associated \emph{Artin monoid} and it can be also defined by the monoid presentation
(see \cite{par})
$$\mathcal{M}(\Gamma )=\big\langle s\in S\,\big|\mathop{\underbrace{sts\dots}}
\limits_{m_{st}\,\text{factors}}=\mathop{\underbrace{tst\dots}}\limits_{m_{st}\,
\text{factors}}\mbox{ for all }s,t\in S\mbox{ satisfying }m_{st}\neq\infty \big\rangle.$$
There are two obvious surjective morphisms:

\begin{picture}(300,45)
\put(160,26){$\mathcal{M}(\Gamma )$}   \put(200,26){$\hookrightarrow$}
\put(192,6){$W(\Gamma)$}.              \multiput(170,23)(3,-3){2}{\vector(1,-1){10}}
\put(220,26){$\mathcal{A}(\Gamma )$}   \multiput(233,23)(-3,-3){2}{\vector(-1,-1){10}}
\end{picture}

Now we recall some basic properties of the (absolute) Garside element of an Artin spherical
monoid $\mathcal{M}(\Gamma )$ (see~\cite{gar},~\cite{bs},~\cite{del},~\cite{mich},~\cite{par},
and also Section 2 for notation). This element is the least common left-multiple of the set of generators
$x_{1},\ldots ,x_{n}$: we have, for any $i=1,\dots,n$, $x_{i}\mid_{L}\Delta $ and if $x_{i}\mid_{L}\omega $ for all
$i$, then $\Delta \mid_{L}\omega $. The element $\Delta (\Gamma )$ is square free (there is no generator
$x_{i}$ such that $x_{i}^{2}\mid\Delta $). In some cases $\Delta (\Gamma )$ itself is a square: for example
$$\Delta (A_{3})=x_{1}(x_{2}x_{1})(x_{3}x_{2}x_{1})=(x_{1}x_{3}x_{2})^{2}; $$
$$\Delta \big(I_{2}(4k)\big)=(\underbrace{x_{1}x_{2}\dots x_{2}}\limits_{2k\,\text{times}})^{2}.$$
Also, there is a bijection $\sigma :\{1,\ldots ,n\}\rightarrow \{1,\ldots,n\}$ such that
$x_{i}\Delta =\Delta x_{\sigma (i)}$. The image of $\Delta(\Gamma )$ in the corresponding
Coxeter group $W(\Gamma )$ is the (unique) longest element of this group and it has order two (see~\cite{bbk} 
and~\cite{DAV}). The length $l(\Gamma )$ of the Garside element $\Delta (\Gamma )$ is equal to the number 
of the reflections (i.e. the conjugates of the generators in the Coxeter group) and it is given by the 
following table (see~\cite{hum} and~\cite{bb}):
\vspace{0.3cm}
\begin{center}
\begin{tabular}{|l|l|l|l|l|l|l|l|l|l|l|l|}
\hline $ \,\,\,\Gamma $ & $ \,\,\,\,A_{n}$ & $B_{n}$ & $\,\,\,\,
D_{n}$ & $E_{6}$ & $E_{7}$ & $\,E_{8}$ & $ F_{4}$ & $G_{2}$ &
$H_{3}$ & $H_{4}$ & $I_{2}(p)$ \\ \hline $l(\Gamma )$ &
${n+1\choose 2}$ & $n^2$  &
$n^2-n$  & 36 & 63 & 120 & 24 & \, 6 & 15 & \,60 & \,\,\,\,\,$p$ \\
\hline
\end{tabular}
\end{center}

\section{Relative Garside Elements}\label{s.rge}

Let us fix the notation. We already used the divisibility relation between two elements of a monoid
$\mathcal{M}$, $\alpha \mid\beta $: this is equivalent to $\beta=\lambda \alpha \rho$ for some elements
$\lambda,\rho\in \mathcal{M}$. If $\lambda=1$, we write $\alpha \mid_L\beta $ and similarly, if
$\rho=1$, we write $\alpha \mid_R\beta $ and we say that $\alpha $ is a left and right divisor of
$\beta $ respectively. The element $\alpha $ in a monoid ${\mathcal{M}}$ generated by $x_1,x_2,\ldots$ is
said to be \emph{rigid} if $\alpha $ can be represented in a unique way as a word in $x_1,x_2,\ldots$.

Suppose we have an inclusion of Coxeter graphs $\Gamma_{n-1}\subset\Gamma_{n}
$ with vertices $\{x_{1},\ldots,x_{n-1}\}$ and $\{x_{1},\ldots,x_{n}\}$, respectively. Using the
definition of Garside elements we have
$$\Delta(\Gamma_{n-1})\mid_{L}\Delta(\Gamma_{n}).$$

\begin{defn}
The relative Garside element $\Delta({\Gamma_{n},\Gamma_{n-1}})$ is defined as a right quotient:
$$\Delta(\Gamma_{n})=\Delta(\Gamma_{n-1})\Delta({\Gamma_{n},\Gamma_{n-1}}).$$
\end{defn}
\noindent If there is no ambiguity concerning the inclusion
$\Gamma_{n-1}\subset\Gamma_{n},$ we will use the simple notation
$\Delta_n=\Delta_{n-1}R_n$. Now we present some properties of the
relative Garside element $R_n$ which characterize this element.

\begin{prop}\label{prop1}
The relative Garside element $R_{n}=\Delta(\Gamma_{n},\Gamma_{n-1})$ satisfies the properties:

a) $R_{n}$ is square free;

b) $x_{i}\mid_{L}\,R_{n}$ if and only if $i=n$;

c) there is a bijection $\sigma: \{1,\ldots,n-1\}\rightarrow \{1,\ldots,n\}\setminus\{\sigma_n(n)\}$ such
that $x_iR_n=R_n x_{\sigma(i)}$;

d) $x_{j}\mid_{R}\,R_{n}$ if and only if $j=\sigma_n(n)$.
\end{prop}

\proof Let $\Delta_{n}=\Delta_{n-1}R_{n}$. Then we have:

a) $R_{n}$ is square free because $\Delta_{n}$ is square free;

b) $x_{i}\,(1\leq i\leq n-1)$ cannot be a left divisor of $R_{n}$: otherwise, $R_{n}=x_{i}U$. But
$x_{i}\mid_{R}\,\Delta_{n-1}$ implies $\Delta_{n-1}=V x_{i}$, therefore we have
$\Delta_{n}=\Delta_{n-1}R_{n}=Vx_{i}^{2}U$, a contradiction.

c) Let $\sigma_{n-1}:\{1,\ldots,n-1\}\rightarrow\{1,\ldots,n-1\}$ and
$\sigma_{n}:\{1,\ldots,n\}\rightarrow\{1,\ldots,n\}$ be the bijections defined by the conjugation with
$\Delta_{n-1}$ and $\Delta_{n}$, respectively: $x_{i}\Delta_{n-1}=\Delta_{n-1}x_{\sigma_{n-1}(i)}$,
$x_{i}\Delta_{n}=\Delta_{n} x_{\sigma_{n}(i)}$. Now, for $i\in\{1,\dots,n-1\}$, we have:
$$\Delta_{n-1}x_{i}R_n=x_{\sigma^{-1}_{n-1}(i)}\Delta_{n-1}R_n=x_{\sigma^{-1}_{n-1}(i)}\Delta_{n}=
\Delta_{n}x_{\sigma_{n}\circ\sigma^{-1}_{n-1}(i)}=\Delta_{n-1}R_{n}x_{\sigma_{n}\circ\sigma^{-1}_{n-1}(i)},$$
therefore $x_{i}R_n=R_nx_{\sigma(i)},$ where $\sigma=\sigma_{n}\circ\sigma^{-1}_{n-1}$. The image of
$\sigma$ contains all the elements $1,\ldots,n$, but not $\sigma_{n}(n)=m$.

d) If $x_{j}\neq x_{m}$, then there exists $i\in\{1,\ldots,n-1\}$ such that $j=\sigma(i)$ and
$x_{i}R_n=R_nx_{j}$. Suppose that $x_{j}\mid_{R}R_n$, i.e. $R_n=Sx_{j}$. We obtain a contradiction:
$x_{i}R_n=R_nx_{j}=Sx^{2}_{j}$, because $x_{i}R_n$ is a right divisor of $\Delta_{n}$.
\endproof

\begin{rk}
From this proposition, the first and the last factors of $R_n$ are uniquely defined.
In some cases $R_n$ is completely rigid, for instance $R_3(A_3,A_2)=x_3x_2x_1$, in other cases only
the interior factors can be changed, for instance $R_3(A_3,A_1\times A_1)=x_2x_1x_3x_2=x_2x_3x_1x_2$.
\begin{center}
\begin{picture}(380,32)
\thicklines                  \put(25,17){\line(1,0){30}}
\multiput(23,14)(30,0){2}{$\bullet$} \put(100,17){\line(1,0){60}}
\multiput(23,7)(75,0){2}{$x_{1}$}
\multiput(98,14)(30,0){3}{$\bullet$} \put(0,15){$A_{2}:$}
\multiput(53,7)(75,0){2}{$x_{2}$}    \put(158,7){$x_{3}$}
\put(63,15){$\subset A_3:$}  \put(180,15){$A_{1}\times A_{1}:$}
\multiput(228,7)(80,0){2}{$x_{1}$} \put(310,17){\line(1,0){60}}
\put(270,15){$\subset A_{3}:$} \multiput(258,7)(110,0){2}{$x_{3}$}
\put(338,7){$x_{2}$} \multiput(228,14)(30,0){2}{$\bullet$}
\multiput(308,14)(30,0){3}{$\bullet$}
\end{picture}
\end{center}
\end{rk}
We will use the notation $m=\sigma_n(n)$.

\begin{prop}\label{prop2}
If $U_{n}\in\mathcal{M}(\Gamma_{n})$ satisfies the next three conditions:

a) $x^2_{m}\nmid\,U_n$,

b) $x_{i}\mid_{L}\,U_{n}$ if and only if $i=n$, and

c) there is a bijection $\tau: \{1,\ldots,n-1\}\rightarrow \{1,\ldots,\widehat{k},\ldots,n\}$ such
that $x_iU_n=U_n x_{\tau(i)}$,

\noindent then $U_n=R_n$ (and also $k=m$, $\tau=\sigma$).
\end{prop}

\proof Let us define $D=\Delta_{n-1}U_n$. Because of the relations $x_i\mid_L\Delta_{n-1},i=1,\ldots,n-1$ and also
from $x_n\mid_LD$, a consequence of b) and c), we have $x_i\mid_{L}\Delta_{n-1}U_n$ for any $i=1,\ldots,n,$ therefore
$\Delta_{n}\mid_LD$. This implies that $\Delta_{n-1}R_n\mid_L\Delta_{n-1}U_n$, hence $R_n\mid_LU_n$. If $R_n\neq U_n$,
then $U_n=R_nx_jW$ for some $j$. If $j\neq m,$ then $j=\sigma_{n}(i)$ for some $i\in \{1,\ldots,n-1\}$, hence
$U_n=R_nx_{\sigma(i)}W=x_iR_nW$ and this contradicts b). If $j=m$, then $U_{n}=(Sx_{m})x_{m}W$ which contradicts a).
\endproof

\begin{prop}\label{prop3}
If $U_{n}\in\mathcal{M}(\Gamma_{n})$ satisfies the next three conditions:

a) $x_{n}\mid_{L}\,U_{n}$,

b) there is a bijection $\tau:\{1,\ldots,n-1\}\rightarrow \{1,\ldots,\widehat{k},\ldots,n\}$ such
that $x_iU_n=U_n x_{\tau(i)}$, and

c) $U_n$ has minimal length among the words satisfying a) and b),

\noindent then $U_n=R_n$ $($and also $k=m$, $\tau=\sigma)$.
\end{prop}

\proof Let us define $D=\Delta_{n-1}U_n$. As in the previous proof, a) and b) implies $x_i\mid_LD$ for all
$i=1,\ldots,n$. Therefore $\Delta_{n}\mid_LD$, hence $R_n\mid_LU_n$. We have $|U_n|\geq|R_n|$ and
$R_n$ satisfies a) and b), therefore condition c) implies $R_n=U_n$.
\endproof

With the same proof, we have another version of the previous
proposition:
\begin{prop}\label{prop4}
If $U_{n}\in\mathcal{M}(\Gamma_{n})$ satisfies the next three conditions:

a) $x_{n}\mid_{L}\,U_{n}$,

b) there is a bijection $\tau:\{1,\ldots,n-1\}\rightarrow \{1,\ldots,\widehat{k},\ldots,n\}$ such
that $x_iU_n=U_n x_{\tau(i)}$, and

c) the length of $U_n$ is the expected length
$l(\Gamma_{n})-l(\Gamma_{n-1})$ (from the table in Section 1),

\noindent then $U_n=R_n$ $($and also $k=m$, $\tau=\sigma)$.
\end{prop}

\section{The Relative Garside and Garside Elements for the infinite series}
\label{compRGaG}

After the computation of various relative Garside elements, we will describe the Garside
elements associated to the classical list of the connected Coxeter diagrams: the infinite series in this section, 
the exceptional ceses in the next section. Using the
following obvious result, these will give formulae for the Garside elements of all Artin
monoids of spherical type.

\begin{lem}
If the graph $\Gamma$ is the disjoint union of $\Gamma_1$ and
$\Gamma_2$, then $\mathcal{M}(\Gamma)\cong
\mathcal{M}(\Gamma_1)\times\mathcal{M}(\Gamma_2)$ and
$$\Delta(\Gamma)=\Delta(\Gamma_1)\Delta(\Gamma_2)=\Delta(\Gamma_2)\Delta(\Gamma_1).$$
\end{lem}

\noindent ${\bf A_n}$ \textbf{series.} We start with $A_1$:
\begin{picture}(25,15)
\put(10,1){$\bullet$} \put(10,-5){$x_{1}$}
\end{picture} and its Garside element $\Delta(A_{1})=\Delta_{2}=x_1$.

\bigskip

\noindent Construction of $\Delta(A_{n},A_{n-1})$: The element $M_{n}=x_{n}x_{n-1}\ldots x_{1}$ satisfies
the conditions of Proposition \ref{prop2}:

a) $M_{n}$ is square free: obvious because $M_{n}$ is rigid;

b) $x_{i}\mid_{L}M_{n}$ if and only if $i=n$ for the same reason;

c) $x_{i}M_{n}=M_{n}x_{i+1}$, therefore $\sigma(i)=i+1$, where $i=1,\ldots,n-1$, and $m=1$.

\begin{cor}(classical Garside element)\label{c1}
$$\Delta(A_n)=\Delta_{n+1}=x_{1}(x_{2}x_{1})\ldots(x_{n-1}\ldots x_{1})(x_{n}\ldots x_{1}).$$
\end{cor}

\noindent ${\bf B_n}$ \textbf{series.} We start with $B_2$:
\begin{picture}(55,22)
\multiput(10,1)(30,0){2}{$\bullet$} \put(10,-5){$x_{1}$}
\put(40,-5){$x_{2}$} \put(12,4){\line(1,0){30}} \put(25,8){$4$}
\end{picture} and its Garside element $\Delta(B_2)=x_1(x_2x_1x_2)$, see the construction $\Delta(I_{2}(p),A_{1})$.

\noindent Construction of $\Delta(B_{n},B_{n-1})$: The element $N_{n}=x_{n}x_{n-1}\ldots x_{2}x_{1}x_{2}\dots x_{n}$
satisfies the conditions of Proposition \ref{prop2}:

a) $N_{n}$ is square free: obvious, because $N_{n}$ is rigid in $\mathcal{M}(B_n)$;

b) $x_{i}\mid_{L}N_{n}$ if and only if $i=n$ for the same reason;

c) $x_{i}N_{n}=N_{n}x_{i}$, therefore $\sigma(i)=i$, where $i=1,\ldots,n-1$, and $m=n$.

\begin{cor}\label{c2}
$$\Delta(B_{n})=x_{1}(x_{2}x_{1}x_{2})(x_{3}x_{2}x_{1}x_{2}x_{3})\dots
(x_{n}x_{n-1}\dots x_{2}x_{1}x_{2}\dots x_{n-1}x_{n}).$$
\end{cor}

\noindent ${\bf D_n}$ \textbf{series.} We start with $D_3$:
\begin{picture}(55,22)
\multiput(30,-9)(0,20){2}{$\bullet$} \put(8,-5){$x_{1}$}
\put(39,-8){$x_{3}$}                 \put(39,13){$x_{2}$}
\put(12,4){\line(2,1){20}}           \put(12,4){\line(2,-1){20}}
\put(10,1){$\bullet$}
\end{picture} and its Garside element $\Delta(D_3)=x_1(x_2x_1)(x_3x_1x_2)$ (this is
$\Delta(A_3,A_2)=x_3x_2x_1$ with a change of notation).

\noindent Construction of $\Delta(D_{n},A_{n-1})$: Consider the
product
$$ P_{n}=(x_{n}x_{n-2}x_{n-3}\ldots x_{1})(x_{n-1}x_{n-2}\ldots x_{2})
(x_{n}x_{n-2}x_{n-3}\ldots x_{3})(x_{n-1}x_{n-2}\ldots
x_{4})\ldots .$$ We prove by induction that $P_n=R_n$. The initial
step is given by the previous formula:
$\Delta(D_3,A_2)=(x_3x_1)x_2$. We check the conditions in
Proposition \ref{prop2} in the case of $n$ even, when the product
ends with the generator $x_n$:
$$ P_n=(x_{n}x_{n-2}x_{n-3}\ldots x_1) (x_{n-1}x_{n-2}x_{n-3}\ldots x_2)
(x_{n}x_{n-2}x_{n-3}\ldots x_{3})\ldots(x_{n-1}x_{n-2})x_n. $$

a) First we prove that $P_n$ is square free and also that
$x_{n}x_{n-1}$ is not a divisor of $P_n$.

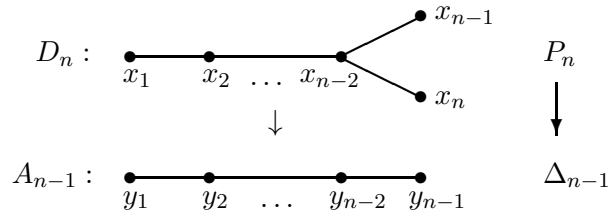
\begin{center}
\begin{picture}(280,110)
\thicklines \put(50,73){\line(1,0){80}}
\put(130,72.5){\line(2,1){30}} \put(130.5,72.5){\line(2,-1){30}}
\put(48,70){$\bullet$}         \put(78,70){$\bullet$}
\put(128,70){$\bullet$}        \put(158,85.5){$\bullet$}
\put(158,55){$\bullet$}        \put(48,63){$x_{1}$}
\put(78,63){$x_{2}$}           \put(115,63){$x_{n-2}$}
\put(166,55){$x_{n}$}          \put(166,87){$x_{n-1}$}
\put(96,62){$\dots$}           \put(15,71){$D_{n}:$}
\put(207,71){$P_n$}            \put(207,25){$\Delta_{n-1}$}
\put(50,27){\line(1,0){110}}   \put(103,45){$\downarrow$}
\put(212,63){\vector(0,-1){20}}\put(48,24){$\bullet$}
\put(78,24){$\bullet$}         \put(158,24){$\bullet$}
\put(128,24){$\bullet$}        \put(48,17){$y_{1}$}
\put(78,17){$y_{2}$}           \put(126,17){$y_{n-2}$}
\put(156,17){$y_{n-1}$}        \put(100,16){$\dots$}
\put(6,25){$A_{n-1}:$}
\end{picture}
\end{center}

\noindent Define a projection $pr:\mathcal{M}(D_n)\rightarrow
\mathcal{M}(A_{n-1})$ given by the diagram:
$pr(x_i)=y_{\min(i,n-1)}$. The image of $P_n$ under the map $pr$
is $\Delta_{n-1}$ which is square free: if $x^{2}_{j}\mid P_n$,
then we have $pr(x^{2}_{j})=y^{2}_{h}\mid\Delta_{n-1}$ ($h=j$ if
$j\leq n-1$ and $h=n-1$ if $j=n$), a contradiction. Therefore
$P_n$ is square free. In the same way, if $x_{n}x_{n-1}\mid P_n$,
then $pr(x_{n}x_{n-1})=y^{2}_{n-1}\mid\Delta_{n-1}$, a
contradiction.

b) Obviously $x_n\mid_LP_{n}$. Let us suppose that
$x_{i}\mid_{L}P_n$ for some $i\leq n-1$. By induction we know that
$R_{n-1}$ is given by the next formula and we will analyze the
following three cases: \newline
1) $i=n-1$, 2) $i=n-2$, and 3) $i\leq n-3$:
$$R_{n-1}=(x_{n-1}x_{n-3}x_{n-4}\ldots x_{1})(x_{n-2}x_{n-3}\ldots x_{2})
(x_{n-1}x_{n-3}x_{n-4}\ldots x_{3})\ldots x_{n-2}. $$ Using the
inclusion morphism $in:\mathcal{M}(D_{n-1})\rightarrow
\mathcal{M}(D_{n})$, $z_1\mapsto x_2$, $z_2\mapsto x_3,\ldots,$
$z_{n-3}\mapsto x_{n-2}$, $z_{n-2}\mapsto x_n$, $z_{n-1}\mapsto
x_{n-1}$, given by the inclusion of graphs:
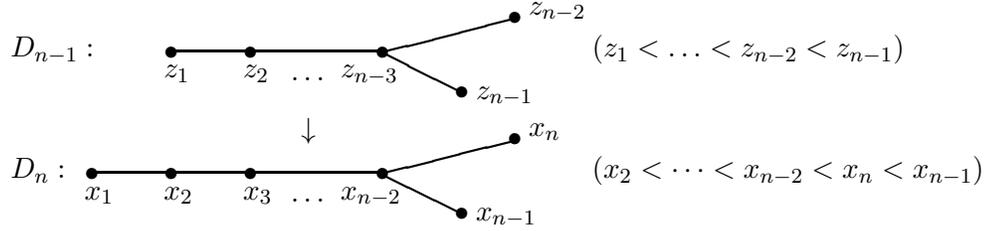
\begin{center}
\begin{picture}(320,60)
\thicklines                    \put(50,27){\line(1,0){80}}
\put(130,26.5){\line(4,1){50}} \put(130.5,26.5){\line(2,-1){30}}
\put(48,24){$\bullet$}         \put(78,24){$\bullet$}
\put(128,24){$\bullet$}        \put(178,37){$\bullet$}
\put(158,9){$\bullet$}         \put(48,17){$z_{1}$}
\put(78,17){$z_{2}$}           \put(115,17){$z_{n-3}$}
\put(166,9){$z_{n-1}$}         \put(186,41){$z_{n-2}$}
\put(96,16){$\dots$} \put(-10,25){$D_{n-1}:$}
\put(210,25){$(z_{1}<\ldots<z_{n-2}<z_{n-1})$}
\end{picture}
\begin{picture}(320,45)
\thicklines \put(20,27){\line(1,0){110}}
\put(130,26.5){\line(4,1){50}} \put(130.5,26.5){\line(2,-1){30}}
\put(48,24){$\bullet$}         \put(78,24){$\bullet$}
\put(128,24){$\bullet$}        \put(178,37){$\bullet$}
\put(158,9){$\bullet$}         \put(48,17){$x_{2}$}
\put(78,17){$x_{3}$}           \put(115,17){$x_{n-2}$}
\put(166,9){$x_{n-1}$}         \put(186,41){$x_{n}$}
\put(96,16){$\dots$}           \put(-10,25){$D_{n}:$}
\put(100,40){$\downarrow$}     \put(18,24){$\bullet$}
\put(18,17){$x_{1}$}
\put(210,25){$(x_{2}<\dots<x_{n-2}<x_{n}<x_{n-1})$}
\end{picture}
\end{center}

\noindent the image of
$$ R_{n-1}=(z_{n-1}z_{n-3}z_{n-4}\dots z_1)(z_{n-2}z_{n-3}\dots z_2)
(z_{n-1}z_{n-3}\dots z_3)\dots(z_{n-1}z_{n-3})z_{n-2} $$
is
$in(R_{n-1})=R'_{n-1}$, given by
$$ R'_{n-1}=(x_{n-1}x_{n-2}x_{n-3}\dots x_{2})(x_{n}x_{n-2}x_{n-3}\dots x_{3})
(x_{n-1}x_{n-2}\dots x_{4})\dots(x_{n-1}x_{n-2})x_{n} $$ and we
have $P_n=(x_{n}x_{n-2}x_{n-3}\dots x_{1})R'_{n-1}$. If
$x_{i}\mid_{L}R'_{n-1}$, then $i=n-1$ (by induction this is true
for $i\in\{2,3,\dots,n,n-1\}$ and also $x_{1}\mid_{L}R'_{n-1}$ is
impossible because $R'_{n-1}$ does not contain $x_{1}$). In the
case 1), $x_{n-1}\mid_{L}P_n$ implies $x_{n-1}x_{n}\alpha_{1}=P_n$
(by Garside Lemma, see Section 5) and $pr(P_{n})$ contains
$y^{2}_{n-1}$, a contradiction. In the case 2),
$x_{n-2}\mid_{L}P_{n}$, and Garside Lemma implies that
$x_{n}x_{n-2}x_{n}\alpha_{2}=x_{n}x_{n-2}x_{n-3}\dots
x_{1}R_{n-1}$, hence $x_{n}\alpha_{2}=x_{n-3}x_{n-4}\dots
x_{1}R_{n-1}$ and $x_{n}x_{n-3}\dots x_{1}\alpha_{3}=x_{n-3}\dots
x_{1}R_{n-1}$, and this gives a contradiction:
$x_{n}\alpha_{3}=R_{n-1}$. In the last case, 3), if
$x_{i}\mid_{L}P_{n}\,(i=1,\dots,n-3)$, then
$x_{n}x_{i}\beta_{1}=x_{n}x_{n-2}\dots x_{i+1}x_{i}\dots
x_{1}R_{n-1}$, and using Garside Lemma we have
$x_{i}\beta_{1}=x_{n-2}\dots x_{i+1}x_{i}\dots x_{1}R_{n-1}$, and
also $x_{i}x_{n-2}\ldots x_{i+2}\beta_{2}=x_{n-2}\ldots
x_{i+2}x_{i+1}x_{i}\ldots x_{1}R_{n-1}$. We obtain
$x_{i}\beta_{2}=x_{i+1}x_{i}\ldots x_{1}R_{n-1}$, next
$x_{i+1}x_{i}x_{i+1}\beta_{3}=x_{i+1}x_{i}x_{i-1}\ldots
x_{1}R_{n-1}$ and $x_{i+1}x_{i-1}\ldots
x_{1}\beta_{4}=x_{i-1}\ldots x_{1}R_{n-1}$; this gives another
contradiction: $x_{i+1}\mid_{L}R_{n-1}$.

c) We have $x_{i}P_{n}=P_{n}x_{n-i}$, hence $\sigma(i)=n-i$.
Therefore $P_{n}=R_{n}$.

Similarly one can check the conditions a), b), c) of
Proposition~\ref{prop2} for $n$ odd:
$$R_{n}=(x_{n}x_{n-2}x_{n-3}\dots x_{1})(x_{n-1}x_{n-2}\dots x_{2})
(x_{n}x_{n-2}x_{n-3}\dots x_{3})\dots(x_{n}x_{n-2})x_{n-1}. $$

\begin{cor}\label{c3}
$$ \Delta(D_{n})=x_{1}(x_{2}x_{1})\dots(x_{n-1}\dots x_{1})
(x_{n}x_{n-2}x_{n-3}\dots x_{1})(x_{n-1}x_{n-2}\dots
x_{2})(x_nx_{n-2}\dots x_3)\dots \, . $$
\end{cor}

\noindent ${\bf I_2(p)}$ \textbf{series,} ${\bf B_2}$ \textbf{and} ${\bf G_2.}$
Construction of $\Delta(I_{2}(p),A_{1})$: Let us define
$Q_{2}(p)=x_{2}x_{1}x_{2}x_{1}x_{2}\dots $,  ($p-1$ factors). This
element satisfies the conditions in Proposition \ref{prop3}:

a) clearly $x_{2}\mid_LQ_{2}(p)$;

b) we have $x_{1}Q_{2}(2p+1)=Q_{2}(2p+1)x_{2}$ ($m=1$) and
$x_{1}Q_{2}(2p)=Q_{2}(2p)x_{1}$ ($m=2$);

c) any relation $x_1V=Vx_{\sigma(1)}$ should involve the unique
defining relation $ x_2x_1x_2\dots\,(p\mbox{ factors})$
$=x_1x_2x_1\dots\,(p\mbox{ factors})$, so the length of $V$ is
greater than or equal to $p-1$ and $Q_2(p)$ has minimal length
among the words satisfying a) and b).

\begin{cor}\label{c4}
$$ \begin{array}{lll}
  \Delta\big(G_{2}\big)    & = & x_{1}x_{2}x_{1}x_{2}x_{1}x_{2}, \\
  \Delta\big(I_{2}(p)\big) & = & x_{1}x_{2}x_{1}x_{2}\dots\,(p \mbox{ factors}).
\end{array}  $$
\end{cor}

\section{The relative Garside elements for the exceptional series}

As a consequence of the results in Section 3 we obtained a new proof for the lengths
of the Garside elements corresponding to the infinite series $A_*$, $B_*$, $D_*$ and $I(*)$ (including
$G_2$). In this section we will use the lengths of the Garside elements to find the relative Garside
elements and also the Garside elements of the monoids corresponding to the exceptional Coxeter graphs
$E_6,E_7,E_8,F_4,H_3$ and $H_4$.

\noindent ${\bf F_4}$ \textbf{case.} We consider the inclusion:

\begin{center}
\begin{picture}(390,32)
\thicklines                         \put(100,16){\line(1,0){60}}
\multiput(98,13)(30,0){3}{$\bullet$}\put(55,13){$B_{3}:$}
\multiput(143,20)(120,0){2}{$4$}    \put(220,16){\line(1,0){90}}
\multiput(98,6)(120,0){2}{$x_{1}$}  \put(308,6){$x_{4}$}
\multiput(128,6)(120,0){2}{$x_{2}$}
\multiput(158,6)(120,0){2}{$x_{3}$} \put(180,13){$\subset F_{4}:$}
\multiput(218,13)(30,0){4}{$\bullet$}
\end{picture}
\end{center}

\noindent Construction of $\Delta(F_4,B_3)$: Let us define
$T_3=x_3x_2x_1x_3x_2x_3$; the Garside element of $B_3$ (with the
change of marking $x_1\leftrightarrow x_3$) is
$\Delta(B_3)=x_1x_2x_1T_3$. Defining $R_4=x_4T_3x_4T_3x_4$, we
have  $x_iR_4=R_4x_i$, $i=1,2,3$, and $l(R_4)=15$; from the table
in Section 1 we have $l(F_4)-l(B_3)=24-9=15$, hence, using
Proposition \ref{prop4}, we obtain
$$ \Delta(F_4,B_3)=x_4T_3x_4T_3x_4. $$

\begin{cor}\label{c5}
$$\Delta(F_4)=x_1x_2x_1(x_3x_2x_1x_3x_2x_3)x_4(x_3x_2x_1x_3x_2x_3)x_4(x_3x_2x_1x_3x_2x_3)x_4. $$
\end{cor}
\newpage
\noindent ${\bf H_{n=3,4}}$ \textbf{series.} We consider the
inclusions:

\begin{center}
\begin{picture}(340,32)
\thicklines                          \put(20,16){\line(1,0){30}}
\multiput(18,13)(30,0){2}{$\bullet$} \put(-20,13){$I_2(5):$}
\multiput(33,20)(80,0){2}{$5$}       \put(100,16){\line(1,0){60}}
\multiput(98,13)(30,0){3}{$\bullet$} \put(63,13){$\subset H_3:$}
\multiput(18,6)(80,0){2}{$x_{1}$}    \put(218,6){$x_{1}$}
\multiput(48,6)(80,0){2}{$x_{2}$}    \put(248,6){$x_{2}$}
\multiput(158,6)(120,0){2}{$x_{3}$}  \put(308,6){$x_{4}$}
\multiput(218,13)(30,0){4}{$\bullet$}\put(180,13){$\subset H_4:$}
\put(233,20){$5$}                    \put(220,16){\line(1,0){90}}
\end{picture}
\end{center}

Construction of $\Delta(H_{3},I_{2}(5))$ and
$\Delta(H_{4},H_{3})$: The element
$$ S_3=(x_3x_2x_1x_2x_1)(x_3x_2x_1x_2)x_3 $$
satisfies the commutation rules $x_1S_3=S_3x_2$, $x_2S_3=S_3x_1$,
and its length is $10$. From the length table we find that
$l(\Delta(H_3,I_2(5)))=15-5=10$, so
$$ \Delta(H_3,I_2(5))=S_3. $$
Similarly, the element
$$ S_4=x_4S_3x_4S_3x_4S_3x_4S_3x_4 $$
verifies $x_iS_4=S_4x_i$, $i=1,2,3$, and it has the expected length:
$l(S_4)=45=60-15=l(H_4)-l(H_3)$, therefore
$$ \Delta(H_4,H_3)=S_4. $$

\begin{cor}\label{c6}
$$\begin{array}{lll}
 \Delta(H_3) & =  & x_1x_2x_1x_2x_1S_3=x_1x_2x_1x_2x_1\cdot (x_3x_2x_1x_2x_1)(x_3x_2x_1x_2)x_3,\\
 \Delta(H_4) & =  & x_1x_2x_1x_2x_1S_3S_4=x_1x_2x_1x_2x_1S_3x_4S_3x_4S_3x_4S_3x_4S_3x_4.
\end{array} $$
\end{cor}

\noindent ${\bf E_{n=6,7,8}}$ \textbf{series.} We consider the
inclusions:

\begin{center}
\begin{picture}(360,130)
\thicklines                          \put(-10,40){$E_6:$}
\multiput(20,45)(0,65){2}{\multiput(0,0)(170,0){2}{
\multiput(-2,0)(30,0){4}{$\bullet$}  \put(0,3){\line(1,0){90}}
\put(60,1){\line(0,-1){25}}          \put(-5,8){$x_{1}$}
\put(25,8){$x_{2}$}                  \put(55,8){$x_{3}$}
\put(85,8){$x_{5}$}                  \put(67,-30){$x_{4}$}
\put(58,-30){$\bullet$}              }}
\put(110,48){\line(1,0){30}}         \put(138,45){$\bullet$}
\put(138,53){$x_6$}                  \put(280,48){\line(1,0){60}}
\multiput(308,53)(0,65){2}{$x_{6}$}  \put(368,118){$x_{8}$}
\multiput(338,53)(0,65){2}{$x_{7}$}  \put(280,113){\line(1,0){90}}
\multiput(308,45)(30,0){2}{$\bullet$}
\multiput(308,110)(30,0){3}{$\bullet$} \put(-10,105){$D_5:$}
\put(-10,75){$\bigcap$}              \put(150,40){$\subset E_7:$}
\put(163,75){$\bigcup$}              \put(160,105){$E_8:$}
\end{picture}
\end{center}

\noindent Construction of $\Delta(E_6,D_5)$, $\Delta(E_7,E_6)$,
and $\Delta(E_8,E_7)$: Let us define the element
$$ V_6=x_6\Delta(D_5,D_4)x_6x_5x_3x_2x_1=(x_6x_5x_3x_2x_1x_4x_3x_2x_5x_3x_4)(x_6x_5x_3x_2x_1). $$
This verifies the commutation relations:
$$ x_1V_6=V_6x_6,x_2V_6=V_6x_5, x_3V_6=V_6x_3,x_4V_6=V_6x_2 \mbox{ and }x_5V_6=V_6x_4. $$
Define also the elements
$$ V_7=x_7V_6(x_7x_6x_5x_3x_2x_4x_3x_5x_6)x_7 \mbox{ and } V_8=x_8V_7x_8V_7x_8. $$
These elements verify the commutation relations:
$$ \begin{array}{cc}
    x_iV_7= \left\{
    \begin{array}{ll}
    V_7x_{7-i}, & i=1,2,5,6, \\
    V_7x_i, & i=3,4,
   \end{array}
   \right.                      &
\mbox{ and }\quad x_iV_8=V_8x_i,\,\,\, i=1,\dots,7.
\end{array}   $$
Counting the lengths we obtain
$$ \begin{array}{l}
   l(V_6)=16=36-20=l(E_6)-l(D_5), \\
   l(V_7)=27=63-36=l(E_7)-l(E_6), \\
   l(V_8)=57=120-63=l(E_8)-l(E_7).
\end{array} $$
From Proposition \ref{prop4} we obtain the relative Garside
elements:
$$ \Delta(E_6,D_5)=V_6,\Delta(E_7,E_6)=V_7,\Delta(E_8,E_7)=V_8. $$
\begin{cor}\label{c7}
$$\begin{array}{lll}
 \Delta(E_6) & =  & \Delta(D_5)V_6=\Delta(A_4)\Delta(D_5,A_4)V_6,         \\
 \Delta(E_7) & =  & \Delta(E_6)V_7=\Delta(A_4)\Delta(D_5,A_4)V_6V_7,      \\
 \Delta(E_8) & =  & \Delta(E_7)V_8=\Delta(A_4)\Delta(D_5,A_4)V_6V_7V_8.
\end{array} $$
\end{cor}

\section{Garside Lemma and few computations}

The next lemma was proved by Garside for the braid monoid (or
$A_n$ series), see~\cite{gar}, and generalized for an arbitrary
Artin monoid by Brieskorn and Saito, see~\cite{bs}:
\begin{lem}\label{Gtype} \textbf{(Garside Lemma)}
Let $W$ be an element in the Artin monoid $\mathcal{M}$ such that
$x_i\mid_LW$  and $x_j\mid_LW$ ($i\neq j$). Then there is an
element $Z\in \mathcal{M}$ such that
$$ W=(\mathop{\underbrace{x_ix_jx_ix_j\dots}} \limits_{m_{ij}\,\,\text{times}})Z=
(\mathop{\underbrace{x_jx_ix_jx_i\dots}}
\limits_{m_{ij}\,\,\text{times}})Z. $$
\end{lem}

Now we give the details for the proof of two commutation relations
described in Section 4. First a short computation:

\begin{lem} In $F_4$ we have $x_1R_4=R_4x_1$. \end{lem}
\proof
The factors which are transformed under Coxeter relations are written in bold characters:
$$ \begin{array}{l}
   {\bf x_1}\cdot x_4(x_3x_2x_1x_3x_2x_3)x_4(x_3x_2x_1x_3x_2x_3)x_4=\\
   =x_4(x_3{\bf x_1}x_2x_1x_3x_2x_3)x_4(x_3x_2x_1x_3x_2x_3)x_4=\\
   =x_4(x_3x_2x_1{\bf x_2}x_3x_2x_3)x_4(x_3x_2x_1x_3x_2x_3)x_4=\\
   =x_4(x_3x_2x_1x_3x_2x_3){\bf x_2}x_4(x_3x_2x_1x_3x_2x_3)x_4=\\
   =x_4(x_3x_2x_1x_3x_2x_3)x_4({\bf x_2}x_3x_2x_3x_1x_2x_3)x_4=\\
   =x_4(x_3x_2x_1x_3x_2x_3)x_4(x_3x_2x_3{\bf x_2}x_1x_2x_3)x_4=\\
   =x_4(x_3x_2x_1x_3x_2x_3)x_4(x_3x_2x_3x_1x_2{\bf x_1}x_3)x_4=\\
   =x_4(x_3x_2x_1x_3x_2x_3)x_4(x_3x_2x_1x_3x_2x_3)x_4\cdot{\bf x_1}.
\end{array} $$
\endproof
And now a long computation:

\begin{lem} In $E_8$ we have $x_7V_8=V_8x_7$. \end{lem}
\proof
In $E_8$ we have the following sequence of equalities:
$$ \begin{array}{lll}
  \alpha  & \equiv & x_3x_2x_4x_3x_2{\bf x_4}=x_3x_2{\bf x_3}x_4x_3x_2={\bf x_2}x_3x_2x_4x_3x_2\equiv\\
          & \equiv & x_2x_3{\bf x_2}x_4x_3x_2=x_2x_3x_4x_3x_2{\bf x_3}=x_2{\bf x_3}x_4x_3x_2x_3= \\
          & =      & x_2x_4x_3x_2{\bf x_4}x_3\equiv\beta,
\end{array} $$
and from the equality $\alpha=\beta$ we get
$$ \begin{array}{lll}
  \gamma  & \equiv & x_5(x_3x_2x_4x_3x_5)({\bf x_3}x_2x_1x_4x_3x_2)=x_5(x_3{\bf x_5}x_2x_4x_3x_5)(x_2x_1x_4x_3x_2)=\\
          & =      & {\bf x_3}x_5(x_3x_2x_4x_3{\bf x_5})(x_2x_1x_4x_3x_2)=x_3x_5(x_3x_2x_4x_3)(x_2x_1x_4{\bf x_5}x_3x_2)\equiv\\
          & \equiv & x_3x_5(x_3x_2x_4x_3)(x_2{\bf x_1}x_4x_5x_3x_2)=x_3x_5(x_3x_2x_4x_3)(x_2x_4{\bf x_1}x_5x_3x_2)\equiv\\
          & \equiv & x_3x_5\alpha{\bf x_1}x_5x_3x_2=x_3x_5\beta{\bf x_1}x_5x_3x_2\equiv x_3x_5(x_2x_4x_3x_2x_4x_3{\bf x_1}x_5x_3x_2)=\\
          & =      & x_3x_5(x_2x_4x_3x_2{\bf x_1}x_4x_3x_5{\bf x_3}x_2)=x_3{\bf x_5}(x_2x_4x_3{\bf x_5})(x_2x_1x_4x_3x_5x_2)=\\
          & =      & x_3(x_2x_4{\bf x_5}x_3x_5)(x_2x_1x_4x_3x_5{\bf x_2})=(x_3x_2x_4x_3x_5)({\bf x_3}x_2x_1x_4x_3{\bf x_2})x_5\equiv\delta,
\end{array} $$
and also, from $\gamma=\delta$, we obtain
\newpage
$$ \begin{array}{lll}
    \eta  & \equiv & x_6(x_5x_3x_2x_4x_3x_5x_6)({\bf x_5}x_3x_2x_1x_4x_3x_2x_5x_3x_4)=\\
          & \equiv & x_6(x_5{\bf x_6}x_3x_2x_4x_3x_5{\bf x_6})(x_3x_2x_1x_4x_3x_2x_5x_3x_4)=\\
          & =      & {\bf x_5}x_6x_5(x_3x_2x_4x_3x_5)(x_3x_2x_1x_4x_3x_2){\bf x_6}x_5x_3x_4\equiv\\
          & \equiv & x_5x_6\gamma x_6x_5x_3x_4=x_5x_6\delta x_6x_5x_3x_4\equiv\\
          & \equiv & x_5x_6(x_3x_2x_4x_3x_5)(x_3x_2x_1x_4x_3x_2)x_5x_6{\bf x_5}x_3x_4=\\
          & =      & x_5x_6(x_3x_2x_4x_3x_5)(x_3x_2x_1x_4x_3x_2{\bf x_6}x_5x_6x_3x_4)=\\
          & =      & x_5{\bf x_6}(x_3x_2x_4x_3x_5{\bf x_6})(x_3x_2x_1x_4x_3x_2x_5{\bf x_6}x_3x_4)=\\
          & =      & (x_5x_3x_2x_4x_3x_5x_6)({\bf x_5}x_3x_2x_1x_4x_3x_2x_5x_3x_4){\bf x_6}\equiv\theta.
\end{array} $$
For the final step we use the next equality
$$ \begin{array}{lll}
  \lambda & \equiv & x_7(x_6x_5x_3x_2x_4x_3x_5x_6x_7)V_6\equiv\\
          & \equiv & x_7(x_6x_5x_3x_2x_4x_3x_5x_6x_7)({\bf x_6}x_5x_3x_2x_1x_4x_3x_2x_5x_3x_4)(x_6x_5x_3x_2x_1)=\\
          & =      & x_7(x_6{\bf x_7}x_5x_3x_2x_4x_3x_5x_6{\bf x_7})(x_5x_3x_2x_1x_4x_3x_2x_5x_3x_4)(x_6x_5x_3x_2x_1)=\\
          & =      & {\bf x_6}x_7x_6(x_5x_3x_2x_4x_3x_5x_6)(x_5x_3x_2x_1x_4x_3x_2x_5x_3x_4)({\bf x_7}x_6x_5x_3x_2x_1)\equiv\\
          & \equiv & x_6x_7\eta x_7x_6x_5x_3x_2x_1=x_6x_7\theta x_7x_6x_5x_3x_2x_1\equiv\\
          & \equiv & x_6{\bf x_7}(x_5x_3x_2x_4x_3x_5x_6)(x_5x_3x_2x_1x_4x_3x_2x_5x_3x_4)x_6x_7{\bf x_6}x_5x_3x_2x_1=\\
          & =      & (x_6x_5{\bf x_7}x_3x_2x_4x_3x_5x_6{\bf x_7})(x_5x_3x_2x_1x_4x_3x_2x_5x_3x_4)(x_6{\bf x_7}x_5x_3x_2x_1)=\\
          & =      & (x_6x_5x_3x_2x_4x_3x_5x_6x_7)({\bf x_6}x_5x_3x_2x_1x_4x_3x_2x_5x_3x_4)(x_6x_5x_3x_2x_1){\bf x_7}\equiv\\
          & \equiv & (x_6x_5x_3x_2x_4x_3x_5x_6x_7)V_6x_7\equiv\mu,
\end{array} $$
and we find
$$ \begin{array}{lll}
   x_7V_8 & \equiv & x_7(x_8V_7x_8V_7x_8)\equiv \\
          & \equiv & {\bf x_7}[x_8x_7V_6(x_7x_6x_5x_3x_2x_4x_3x_5x_6){\bf x_7}x_8x_7V_6(x_7x_6x_5x_3x_2x_4x_3x_5x_6)x_7x_8]=\\
          & =      & x_8x_7V_6({\bf x_8}x_7x_6x_5x_3x_2x_4x_3x_5x_6){\bf x_8}x_7V_6({\bf x_8}x_7x_6x_5x_3x_2x_4x_3x_5x_6)x_7x_8=\\
          & =      & x_8x_7V_6{\bf x_7}x_8(x_7x_6x_5x_3x_2x_4x_3x_5x_6)x_7V_6)(x_8x_7x_6x_5x_3x_2x_4x_3x_5x_6)x_7x_8\equiv\\
          & \equiv & x_8x_7V_6x_7x_8\lambda(x_8x_7x_6x_5x_3x_2x_4x_3x_5x_6)x_7x_8=\\
          & =      & x_8x_7V_6x_7x_8\mu(x_8x_7x_6x_5x_3x_2x_4x_3x_5x_6)x_7x_8\equiv\\
          & \equiv & x_8x_7V_6x_7x_8(x_6x_5x_3x_2x_4x_3x_5x_6x_7V_6)(x_7x_8{\bf x_7}x_6x_5x_3x_2x_4x_3x_5x_6)x_7x_8=\\
          & =      & x_8x_7V_6(x_7{\bf x_8}x_6x_5x_3x_2x_4x_3x_5x_6)x_7{\bf x_8}V_6(x_7{\bf x_8}x_6x_5x_3x_2x_4x_3x_5x_6)x_7x_8=\\
          & =      & [x_8x_7V_6(x_7x_6x_5x_3x_2x_4x_3x_5x_6)x_7x_8{\bf x_7}V_6(x_7x_6x_5x_3x_2x_4x_3x_5x_6)x_7x_8]{\bf x_7}\equiv\\
          & \equiv &(x_8V_7x_8V_7x_8) x_7\equiv V_8x_7.
\end{array} $$
\endproof



\begin{thebibliography}{99}

\bibitem{art} E. Artin, \emph{Theory of braids}, Ann. Math. \textbf{48} (1947), 101-126.

\bibitem{tez} B. Berceanu, \emph{Artin algebras -- applications in topology} (in Romanian) Ph. Thesis,
   University of Bucharest (1995).

\bibitem{BP} B. Berceanu, S. Parveen, \emph{Braid groups in complex projective spaces}, Advances in 
Geometry \textbf{12} (2012), 269-286.
   
\bibitem{bb} A. Bj\"orner, F. Brenti, \emph{Combinatorics of Coxeter Groups}, Graduate Texts
in Mathematics \textbf{231}, Springer Verlag, 2005.

\bibitem{bok} L. A. Bokut, Y. Fong, W. F. Ke, L. S. Shiao.  \emph{Gr\"{o}bner-Shirshov
bases for braid semigroup}, Advances in Algebra, World Sci. Publ.
(2003), 60-72.

\bibitem{bbk} N. Bourbaki, \emph{Groupes et Alg\`{e}bres de Lie}, Chapitres 4-6,
Elem. Math., Hermann, 1968.

\bibitem{bs} E. Brieskorn, K. Saito, \emph{Artin groups and Coxeter groups},
Invent. Math \textbf{17} (1972), 245-271.

\bibitem{DAV} M. Davis, \emph{The Geometry and Topology of Coxeter Groups}, London Mathematical 
Society Monographs, Princeton Univ. Press, 2008.

\bibitem{del} P. Deligne, \emph{Les immeubles des groupes de tresses g\'{e}n\'{e}ralis\'{e}s},
Invent. Math. \textbf{17} (1972),  273-302.

\bibitem{gar} F. A. Garside, \emph{The braid groups and other groups}, Quart.
J. Math. Oxford, $2^e$ ser. \textbf{20} (1969), 235-254.

\bibitem{hum} J. E. Humphreys, \emph{Reflection Groups and Coxeter Groups}, Cambridge Univ. Press, 1990.

\bibitem{mich} J. Michel, \emph{A note on words in braid monoids}, J.
Algebra \textbf{215} (1999), no. 1, 366-377.

\bibitem{mor} S. Moran, \emph{The Mathematical Theory of Knots and Braids}, North-Holland
Mathematical Studies, Vol. \textbf{80}, Elsevier, Amsterdam, 1983.

\bibitem{par}  L. Paris, \emph{Braid groups and Artin groups}, in "Handbook on Teichm\"{u}ller Theory,
Vol \textbf{II}, EMS Publishing House, Z\"{u}rich, 2008.

\end{thebibliography}
\end{document}